\documentclass[12pt]{article}
\usepackage{amssymb,amsmath,amsthm}
\date{\small \jobname: Mar 13, 2009}

\normalsize

\makeatletter
\newfont{\footsc}{cmcsc10 at 8truept}
\newfont{\footbf}{cmbx10 at 8truept}
\newfont{\footrm}{cmr10 at 10truept}
\renewcommand{\ps@plain}{%
\renewcommand{\@oddfoot}{\footsc file: \jobname.tex\hfil\footrm\thepage}} 
\makeatother
\iftrue \makeatletter
\@ifundefined{every@math@size}{}{%
\addto@hook\every@math@size{\dch@scr@hook} 
\def\dch@scr@adjust{\@ifundefined{dch@sizet\f@size}%
  {\expandafter\dch@set@script\csname dch@sizet\f@size\endcsname}%
  {\csname dch@sizet\f@size\endcsname}} 
\def\dch@set@script#1{\begingroup %
  \frozen@everymath{}%
  \let#1\@empty \let\dch@do@one\relax 
  \dch@set@one\scriptscriptstyle\scriptscriptfont#1\ssf@size 
  \dch@set@one\scriptstyle\scriptfont#1\sf@size 
  \dch@set@one\textstyle\textfont#1\f@size 
  \endgroup #1} %
\def\dch@set@one#1#2#3#4{%
  \@ifundefined{dch@size#4}%
   {\expandafter\xdef\csname dch@size#4\endcsname{%
      \fontdimen13\the#2\tw@=\the\fontdimen13#2\tw@ 
      \fontdimen14\the#2\tw@=\the\fontdimen14#2\tw@ 
      \fontdimen15\the#2\tw@=\the\fontdimen15#2\tw@ 
      \fontdimen16\the#2\tw@=\the\fontdimen16#2\tw@ 
      \fontdimen17\the#2\tw@=\the\fontdimen17#2\tw@}%
  }{\csname dch@size#4\endcsname}%
  \setbox\z@\hbox{$#1H_2$}\@tempdima\dp\z@ 
  \setbox\z@\hbox{$#1H_2^{+\vrule \@height 1em}$}%
   \ifdim\@tempdima<\dp\z@ 
    \advance\@tempdima\dp\z@ \divide\@tempdima\tw@ %
    \@tempdimb\dp\z@ \advance\@tempdimb-\@tempdima %
    \advance\@tempdimb\ht\z@ \advance\@tempdimb-1em %
    \xdef#3{#3\dch@do@one#2{\the\@tempdimb}{\the\@tempdima}}%
  \fi} 
\def\dch@do@one#1#2#3{%
  \fontdimen13#1\tw@#2\relax 
  \fontdimen14#1\tw@\fontdimen13#1\tw@ \fontdimen15#1\tw@\fontdimen13#1\tw@ 
  \fontdimen\sixt@@n#1\tw@#3\fontdimen17#1\tw@\fontdimen\sixt@@n#1\tw@}%
\let\dch@scr@hook\dch@scr@adjust 
\ifx\glb@currsize\f@size \dch@scr@adjust \fi}%
\makeatother \fi

\def\tO{\widetilde{O}}

\def\abs#1{\mathopen|#1\mathclose|}

\def\dfrac#1#2{\lower0.15ex\hbox{\large$\frac{#1}{#2}$}}
\def\tbd{(1+\l)^{-1}\,}
\def\({\bigl(}
\def\){\bigr)}

\newcommand{\lab}[1]{\label{#1}}                

\newtheorem{thm}{Theorem}
\newtheorem{cor}{Corollary}
\newtheorem{lemma}{Lemma}
\newtheorem{conj}{Conjecture}
\newtheorem{rem}{Remark}

\def\sumpp{\sum\nolimits''}
\def\sumsum{\sum\nolimits_{j,k}}
\def\Oh{O}

\def\lcm{{\rm lcm}}

\def\A{{\cal A}}

\def\M{{\cal M}}
\def\P{{\cal P}}
\def\Q{{\cal Q}}
\def\R{{\cal R}}
\def\S{{\cal S}}
\def\U{{\cal U}}

\def\eps{\varepsilon}
\def\l{\lambda}

\def\d{\delta}

\def\t{\theta}
\def\p{\phi}

\def\bix{\boldsymbol{x}}

\def\bit{\boldsymbol{\theta}}
\def\bip{\boldsymbol{\phi}}
\def\bis{\boldsymbol{\sigma}}
\def\bitau{\boldsymbol{\tau}}
\def\bithat{\hat{\bit}}
\def\biphat{\hat{\bip}}
\def\that{\hat\t} \let\th=\that
\def\phat{\hat\p} \let\ph=\phat
\def\zvec{\boldsymbol{z}}

\def\tav{\bar{\bit}} \let\tb=\tav
\def\pav{\bar{\bip}} \let\pb=\pav
\def\xav{\bar{\bix}}

\def\Ntheta{N_{\bit}}
\def\Nphi{N_{\bip}}

\def\mw#1{\hat{#1}}
\def\mwA{\mw{A}}
\def\mwa{\mw{a}}
\def\mwB{\mw{B}}
\def\mwC{\mw{C}}

\def\mwE{\mw{E}}
\def\mwF{\mw{F}}

\def\mwJ{\mw{J}}
\def\mwZ{\mw{Z}}

\def\Chat{\hat C}

\def\Re{\operatorname{Re}}
\def\Im{\operatorname{Im}}

\def\Prob{\operatorname{Prob}}

\def\Reals{{\mathbb{R}}}

\makeatletter
\@addtoreset{equation}{section}

\makeatother

\def\nicebreak{\vskip 0pt plus100pt\penalty-100\vskip 0pt plus-100pt}

\begin{document}

\title{
Asymptotic Enumeration of Integer
Matrices \\ with Constant Row and Column Sums}

\author{
E.~Rodney~Canfield\thanks
 {Research supported by the NSA Mathematical Sciences Program} \\
\small Department of Computer Science\\[-0.8ex]
\small University of Georgia\\[-0.8ex]
\small Athens, GA 30602, USA\\[-0.3ex]
\small\texttt{erc@cs.uga.edu}
\and
Brendan~D.~McKay\thanks
 {Research supported by the Australian Research Council}\\
\small Department of Computer Science\\[-0.8ex]
\small Australian National University\\[-0.8ex]
\small Canberra ACT 0200, Australia\\[-0.3ex]
\small\texttt{bdm@cs.anu.edu.au}
}

\maketitle

\begin{abstract}
 Let $s,t,m,n$ be positive integers such that $sm=tn$.
 Let $M(m,s;n,t)$ be the number of $m\times n$ matrices over
 $\{0,1,2,\dots\}$
 with each row summing to $s$ and each column summing to~$t$.
 Equivalently, $M(m,s;n,t)$ counts 2-way contingency tables
 of order $m\times n$ such that the row marginal sums are all~$s$
 and the column marginal sums are all~$t$.
 A third equivalent description is that $M(m,s;n,t)$ is the number
 of semiregular labelled bipartite
 multigraphs with $m$ vertices of degree $s$ and $n$ vertices
 of degree~$t$.
 When $m=n$ and $s=t$ such matrices are also referred to as
 $n \times n$
 magic or semimagic squares with line sums equal to $t$.
 We prove a precise asymptotic formula
 for $M(m,s;n,t)$ which is valid over a range of $(m,s;n,t)$
 in which
 $m,n\rightarrow\infty$ while remaining approximately equal
 and the average entry is not too small.
 This range includes the case where $m/n$, $n/m$, $s/n$
 and $t/m$ are bounded from below.
\end{abstract}

\section{Introduction}\label{section:intro}
Let $m,s,n,t$ be positive integers such that $ms=nt$.  
Let $M(m,s;n,t)$ be the number of $m\times n$ matrices over
$\{0,1,2,\dots\}$
with each row summing to $s$ and each column summing to~$t$.
Figure~\ref{figure} shows an example.

\begin{figure}[ht]
\centering
  \begin{tabular}{cccc|c}
     5&4&0&11&20\\ 2&8&9&1&20\\ 8&3&6&3&20\\
     \hline 15&15&15&15
  \end{tabular}
  \caption{Example of a table counted by $M(3,20;4,15)$.\label{figure}}
\end{figure}

The matrices counted by $M(m,s;n,t)$ arise frequently
in many areas of mathematics, for example
enumeration of permutations with
respect to descents, symmetric function theory, and statistics.
The last field in particular has an extensive literature in which
such matrices are studied as \textit{contingency tables}
or \textit{frequency tables};
see \cite{gail} as an example of enumerative work,
\cite{Greselin} for a lengthy bibliography,
and \cite{DG} for a survey of applications.
The matrices counted by $M(m,s;n,t)$ have also regularly
been the topic
of papers whose primary focus is algorithmic, both with respect
to the problem of generating contingency tables with prescribed
margins at random,
and with respect to approximate counting.
For random generation see, for example \cite{DKM, morris};
for approximate counting see the recent~\cite{bsy}.
Other important recent studies of both random sampling and approximate
counting include~\cite{CDHL, HJ}.

A historically early study of integer matrices with specified row
and column sums appeared in MacMahon \cite{macmahon}, where the
numbers $M(m,s;n,t)$ are identified as coefficients when certain
functions are expanded in terms of standard bases of symmetric
functions.
In a later paper~\cite{macmahon2}, MacMahon studied ``general
magic squares''---$n\times n$ integer matrices with all row
and column sums equal to a prescribed value~$t$.
Stanley \cite[Chap.\,4]{stanleyA} refers to these as
magic squares with line sums equal to $t$, and Stanley's
\cite[Chap.\,1]{stanley} provides further history of this topic.
For fixed~$n$, $M(n,t;n,t)$ is a polynomial in $t$ known as the
\textit{Ehrhart polynomial\/} of the Birkhoff polytope.
This has been computed in closed form for $n\le 9$~\cite{BP2003}.
(The Birkhoff polytope is the set of all doubly stochastic real
matrices.)
The results in the present paper enable us to give the
first asymptotic formula for the volume of this famous
polytope~\cite{CMComing}.

Our focus in this paper is the asymptotic value
of $M(m,s;n,t)$.
The first significant result on asymptotics
was that of Read~\cite{read},
who obtained the asymptotic behavior for $s=t=3$.  In the four year period
from 1971 to 1974 three published papers \cite{bbk, bender,
everett} gave the asymptotic value for the case $s,t$ bounded.
To our knowledge, this was the extent of
published, rigorously proven formulas until the present work.
In a paper currently in preparation~\cite{GMsparse},
Greenhill and McKay obtain
a formula valid in the sparse range $st=o\((mn)^{1/2}\)$.

Define the {\it density\/} $\l=s/n=t/m$ to be the average entry in
the matrix.
As early as~\cite[p.\,100]{good50}, Good implicitly 
gave the estimate $M(m,s;n,t) \approx G(m,s;n,t)$,
where\goodbreak
$$
 G(m,s;n,t) = 
     \frac{\displaystyle\binom{n+s-1}{s}^{\!m}\binom{m+t-1}{t}^{\!n}}
          {\displaystyle\binom{mn+\l mn-1}{\l mn}}\,.
$$
Good spelt this out explicitly in~\cite{good76} and again with Crook
in~\cite{good}.  In the first paper, Good gave a heuristic argument
for this approximation based on steepest-descent, ``leaving aside finer
points of rigor.'' His calculation also yielded another
approximation~\cite[Eqn.\,B2.22]{good76} that is very similar when
$m\approx n$ but larger otherwise.
We will show in this paper that $G(m,s;n,t)$ is a remarkably accurate
approximation of $M(m,s;n,t)$, being out by a constant factor
over a wide range and perhaps always.
Another important estimate of $M(m,s;n,t)$ was developed by
Diaconis and Efron~\cite{DE}.  It aims for accuracy when the
sum of the matrix entries $\lambda mn$ is very large, and is
quite poor when that condition is not met.

Our main result is as follows.

\nicebreak
\begin{thm}\label{thm:main}
Let $s=s(m,n),t=t(m,n)$ be positive integers satisfying $ms=nt$.  Define
$\l = s/n = t/m$.  Let $a,b>0$ be constants
such that $a+b<\frac12$.  Suppose that $m,n \rightarrow \infty$
in such a way that
\begin{equation} \label{Hyp}
\frac{(1+2\l)^2}{4\l(1+\l)}
   \biggl( 1 + \frac{5m}{6n} + \frac{5n}{6m} \biggr)
\le a \log n.
\end{equation}
Then
\begin{align}
M(m,s;n,t) &=
     G(m,s;n,t) \exp\(\tfrac12 + O(n^{-b})\)\label{mainformula} \\
     &= \frac{\(\l^{-\l}(1+\l)^{1+\l}\)^{mn}}
             {(4\pi A)^{(m+n-1)/2} m^{(n-1)/2} n^{(m-1)/2}}\notag\\
     &{\kern 10mm}\times
             \exp\biggl(\frac12 - \frac{1+2A}{24A}
             		\Bigl(\frac mn + \frac nm\Bigr) + O(n^{-b})\biggl).
		\label{mainformula2}
\end{align}
\end{thm}

\begin{cor}\label{biglambda}
Under the conditions of Theorem~\ref{thm:main}, if in addition
$mn/\l^2\to 0$,
\[
   M(m,s;n,t) = \(\l+\tfrac12\)^{(m-1)(n-1)}
      \frac{(mn)!}{m!^n n!^m} \exp\(\tfrac12+O(mn/\l^2+n^{-b})\).
\]
\end{cor}

\begin{rem}\label{RueDeRemarques}
The following is inspired by a similar observation 
made by Good~\cite{good76}.
The number of $k$-tuples of nonnegative integers $n_i$ satisfying
$n_1+\cdots+n_k=K$ is well known to be
$\binom{K+k-1}{k-1}=\binom{K+k-1}{K}$.  This allows an instructive
interpretation of~\eqref{mainformula} as
$M(m,s;n,t) = NP_1P_2E$, where
\begin{align*}
   N &= \binom{mn+\l mn-1}{\l mn},~
   P_1 = N^{-1}\binom{n+s-1}{s}^{\!m}\negthickspace,~
   P_2 = N^{-1}\binom{m+t-1}{t}^{\!n}\negthickspace,\\[0.6ex]
   E &= E(m,s;n,t) = \exp\(\tfrac12 + O(n^{-b})\).
\end{align*}
Clearly $N$ is the number of tables whose
entries sum to $\l mn = sm = tn$.
In the uniform probability space on these $N$ tables, $P_1$
is the probability of the event that all the row sums are
equal to~$s$,
and $P_2$ is the probability of the event that all the column
sums are equal to~$t$.
The final quantity $E$ is thus a correction to account for
the non-independence of these two events.
\end{rem}

In an earlier paper \cite{CM} we showed how $B(m,s;n,t)$, the number of
such matrices whose entries are taken from $\{0,1\}$, can be computed
exactly for small $m,n$.  That algorithm is easily adaptable for efficient
computation of exact values of $M(m,s;n,t)$, and the values so
computed suggest the following conjecture.

\begin{conj}\label{conjecture}
Consider a 4-tuple of positive integers $m,s,n,t$ such
that $ms=nt$.   Define $\Delta(m,s;n,t)$ by
\begin{equation}\label{conjformula}
\begin{split}
M(m,s;n,t) &= G(m,s;n,t)\,
     \Bigl(\frac{m+1}{m}\Bigr)^{(m-1)/2} 
     \Bigl(\frac{n+1}{n}\Bigr)^{(n-1)/2} \\
  &{\kern 40mm}\times
     \exp\Bigl(-\dfrac12 + \frac{\Delta(m,s;n,t)}{m+n}\Bigr).
\end{split}
\end{equation}
Then\/ $0<\Delta(m,s;n,t)<2$.
\end{conj}
The factor $(1+1/m)^{(m-1)/2}$, which approaches
$e^{1/2}$ as $m\to\infty$, and the similar factor 
$(1+1/n)^{(n-1)/2}$, appear naturally in the analysis of
$M(m,s;n,t)$ when one of $m,n$ goes to~$\infty$ much faster than
the other.  This can be seen
in~\cite[eqn.\,(3.5)]{gail} and will be extended
rigorously in a forthcoming paper~\cite{oblong}.
Conjecture~\ref{conjecture} has been proved in several cases: (a)
for $m=n\le 9$, using the exact values from~\cite{BP2003}; (b)
for sufficiently large $m,n$ when $st=o\((mn)^{1/5}\)$, using the
asymptotics derived in~\cite{GMsparse};
(c) for several thousand values of $(m,s;n,t)$ for $m,n\le 30$.
It has also been established to a high degree of confidence
for many larger sizes using a simulation similar to the one
described in~\cite{CM}, which was a variation on a method
of Chen, Diaconis, Holmes and Liu~\cite{CDHL}.


An interesting point of contrast between the integer and
0-1 cases is suggested by the above.  Let $E_0(m,s;n,t)$ be
the quantity corresponding to $E(m,s;n,t)$ in
Remark~\ref{RueDeRemarques} when the number of binary tables
is decomposed in the same manner (see \cite{CGM}).
It was proved by Ordentlich and Roth~\cite{Ord}
that $E_0(m,s;n,t)\le 1$.
A corollary of Corollary~\ref{conjecture} would be that the
opposite bound $E(m,s;n,t)\ge 1$ holds in the integer case.
(In fact the smallest value we have found for $s,t>0$ is 
$E(2,3;3,2)=\tfrac{539}{450}$.)

Some numerical examples comparing our estimates to other estimates
and the correct values appear in Table~\ref{numeric}.

\begin{table}
\caption{Comparison of six estimates for $M(m,s;n,t)$.}
\label{numeric}
\small
\centering
\def\e#1e#2 {$#1{\times}10^{#2}$}
\begin{tabular}{c|ccc|c}
  $m,s,n,t$ &
  $G(m,s;n,t)$ &
  \cite[(B2.22)]{good76} &
  \cite{DE}\\
  &
  \cite[Cor. 5.1]{GMsparse} & 
  \eqref{mainformula} &
  Conjecture~\ref{conjecture}&Exact\\
 \hline
{\vrule height2.8ex width0pt depth0pt}
3,100,3,100 &
\e1.019e7 &\e1.192e7 &\e1.262e7 \\
 & \e1.022e7 &\e1.680e7 &\e(1.316{\pm}0.217)e7 &13268976 \\[0.7ex]
3,98,49,6 &
\e7.594e67 &\e3.716e68 &
\e1.278e68 \\
 & \e9.630e67 &\e1.252e68 &\e(1.017{\pm}0.020)e68 &\e1.01100e68
\\[0.7ex]
3,99,9,33 &
\e2.116e21 &\e2.788e21 &\e2.864e21 \\
 & \e2.506e21 &\e3.488e21 &\e(2.844{\pm}0.236)e21 &\e2.79207e21
\\[0.7ex]
10,20,10,20 &
\e7.434e58 &\e9.021e58 &\e1.511e59 \\
 & \e1.059e59 &\e1.226e59 &\e(1.119{\pm}0.056)e59 &\e1.09747e59
\\[0.7ex]
18,13,18,13 &
\e5.157e127 &\e6.962e127 &\e5.109e130 \\
 & \e7.850e127 &\e8.502e127 &\e(8.065{\pm}0.224)e127 &
\e7.94500e127 \\[0.7ex]
30,3,30,3 &
\e1.404e92 &\e7.496e92 &\e1.086e138 \\
 & \e2.223e92 &\e2.315e92 &\e(2.242{\pm}0.037)e92 &\e2.22931e92
\end{tabular}
\end{table}

Throughout the paper, the asymptotic notation $O(f(m,n))$ refers
to the passage of $m$ and $n$ to $\infty$.  Generally, the constant
implied by the $O$ may depend on $a$, $b$, and the sufficiently small
and fixed constant $\eps$ introduced later in the definition
(\ref{RDef}) of region $\R$.
We also use a modified notation $\tO(f(m,n))$.  A function
$g(m,n)$ belongs to
this class provided that for some constant $c$, which is independent
of $a$, $b$, and $\eps$, we have
$$
g(m,n) = O(f(m,n)n^{c\eps}).
$$
Under the hypothesis (\ref{Hyp}) of Theorem \ref{thm:main}, we have
$\log m\sim\log n$, $n = O(m\log m)$, $m=O(n\log n)$,
and $\l^{-1}=O(\log n)$.  Consequently, in the newly introduced notation,
$n=\tO(m)$, $m=\tO(n)$, and $\l^{-1}
=\tO(1)$.  In general, if $c_1,c_2,c_3,c_4$
are constants, then $m^{c_1+c_2\eps}n^{c_3+c_4\eps}=\tO(m^{c_1}n^{c_3})
=\tO(n^{c_1+c_3})=\tO(m^{c_1+c_3})$.

\nicebreak
\section{An integral for $M(m,s;n,t)$}\label{section:integral}

We express $M(m,s;n,t)$ as an integral
in $(m{+}n)$-dimensional complex space then estimate its value
by the saddle-point method.

\medskip

It is clear that $M=M(m,s;n,t)$ is the coefficient of
$x_1^s\cdots x_m^s\, y_1^t\cdots y_n^t$ in
$$\prod_{j=1}^m \,\prod_{k=1}^n \, \(1 - x_jy_k\)^{\!-1}.$$
Applying Cauchy's Theorem we have

\begin{equation}\label{cauchy1}
 M = \frac{1}{(2\pi i)^{m+n}} \oint\cdots\oint
   \frac{\prod_{j,k}(1-x_j y_k)^{-1}}{
	x_1^{s+1}\cdots x_m^{s+1} y_1^{t+1}\cdots y_n^{t+1}}
   \, dx_1\cdots dx_m\,dy_1\cdots dy_n,
\end{equation}
where each contour circles the origin once in the anticlockwise direction.

It will suffice to take the contours to be circles; specifically, we will
put $x_j=re^{i\t_j}$ and $y_k=re^{i\p_k}$ for each $j,k$, where
$$ r = \sqrt{\frac{\l}{1+\l}}\,. $$
This gives
\begin{equation}\label{Idef1}
 M = \frac{1}{(2\pi)^{m+n}} \(\l^{-\l}(1+\l)^{1+\l}\)^{mn}\, I(m,n),
\end{equation}
where
\begin{equation}\lab{Idef2}
 I(m,n) = \int_{-\pi}^\pi\cdots \int_{-\pi}^\pi
 \frac{\prod_{j,k} \( 1 - \l(e^{i(\t_j+\p_k)}-1)\)^{-1}}{
    e^{is\sum_j\t_j + it\sum_k\p_k}}
   \, d\bit\, d\bip,
\end{equation}
where $\bit=(\t_1,\dots,\t_m)$ and $\bip=(\p_1,\dots,\p_n)$.
Let $F(\bit,\bip)$ denote the integrand in equation (\ref{Idef2}).

\nicebreak
\section{Evaluating the integral}\label{section:evalint} 
This section follows closely the corresponding section in
\cite{CGM}. 
However, all proofs are intended to be complete, and not require
the reader to refer to the latter work,
except for the following lemma
which we restate here without proof.

\smallskip

\begin{lemma}\label{MW3}
Let $\eps', \eps'', \eps''', \bar\eps,\Delta$ be constants
such that $0<\eps'<\eps''<\eps'''$, $\bar\eps\ge 0$,
and $0<\Delta<1$.  The following is true if $\eps'''$
and $\bar\eps$ are sufficiently small.\endgraf
Let $\mwA=\mwA(N)$ be a real-valued function such that
$\mwA(N)=\Omega(N^{-\eps'})$.
Let $\mwa_j=\mwa_j(N)$, $\mwB_j=\mwB_j(N)$,
$\mwC_{jk}=\mwC_{jk}(N)$, $\mwE_j=\mwE_j(N)$,
$\mwF_{jk}=\mwF_{jk}(N)$ and $\mwJ_j=\mwJ_j(N)$ be
complex-valued functions $(1 \le j,k \le N)$
such that $\mwB_j,\mwC_{jk},
\mwE_j,\mwF_{jk}=O(N^{\bar\eps})$,
$\mwa_j=O(N^{1/2+\bar\eps})$, and $\mwJ_j=O(N^{-1/2+\bar\eps})$,
uniformly over $1\le j,k\le N$.
Suppose that
\begin{align*}
f(\zvec) &= \exp\biggl(
   -\mwA N\sum_{j=1}^N z_j^2
   + \sum_{j=1}^{N} \mwa_j z_j^2 + N \sum_{j=1}^{N} \mwB_j z_j^3
   + \sum_{j,k=1}^N \mwC_{jk}z_j z_k^2
\\
&\kern14mm{}
   + N \sum_{j=1}^N \mwE_j z_j^4
   + \sum_{j,k=1}^N \mwF_{jk} z_j^2 z_k^2
  + \sum_{j=1}^{N} \mwJ_j z_j + \delta(\zvec) \biggr)
\end{align*}
is integrable for $\zvec=(z_1,z_2,\ldots,z_N)\in U_N$
and $\delta(N)=\max_{\zvec\in U_N} \abs{\delta(\zvec)} = o(1)$,
where
\[
U_N = \bigl\{ \zvec \bigm| \abs{z_j} \le N^{-1/2+\hat\eps}
  \text{ for\/ $1\le j\le N$}\bigr\},
\]
where $\hat\eps=\hat\eps(N)$ satisfies
$\eps''\le2\hat\eps\le\eps'''$.
Then, provided the $O(\,)$ term in the following
converges to zero,
\[
 \int_{U_N}f(\zvec)\,d\zvec
  = \biggl(\frac{\pi}{\mwA N}\biggr)^{\!N/2}
   \!\exp\(
     \Theta_1+\Theta_2 + O\((N^{-\Delta}+\delta(N)) \mwZ\)
        \),
\]
where
\begin{align*}
  \Theta_1 &= \frac{1}{2\mwA N}\sum_{j=1}^{N}\mwa_j
     + \frac{1}{4\mwA^2N^2}\sum_{j=1}^{N}\mwa_j^2
     + \frac{15}{16\mwA^3N}\sum_{j=1}^{N}\mwB_j^2
     + \frac{3}{8\mwA^3N^2}\sum_{j,k=1}^{N}\mwB_j\mwC_{jk} \\
  &\quad{} + \frac{1}{16\mwA^3N^3}\sum_{j,k,\ell=1}^N \mwC_{jk}\mwC_{j\ell}
    + \frac{3}{4\mwA^2N}\sum_{j=1}^N \mwE_j
    + \frac{1}{4\mwA^2N^2}\sum_{j,k=1}^N \mwF_{jk}
  \displaybreak[0]\\
  \Theta_2 &= \frac{1}{6\mwA^3N^3}\sum_{j=1}^N \mwa_j^3
 + \frac{3}{2\mwA^3N^2}\sum_{j=1}^N \mwa_j\mwE_j
 + \frac{45}{16\mwA^4N^2}\sum_{j=1}^N \mwa_j\mwB_j^2 \\
 &\quad{} + \frac{1}{4\mwA^3N^3}\sum_{j,k=1}^N (\mwa_j+\mwa_k) \mwF_{jk}
       + \frac{3}{4\mwA^2N}\sum_{j=1}^N \mwB_j\mwJ_j
       + \frac{1}{4\mwA^2N^2}\sum_{j,k=1}^N \mwC_{jk}\mwJ_j \\
 &\quad{} + \frac{1}{16\mwA^4N^4}\sum_{j,k,\ell=1}^N
 			(\mwa_j+2\mwa_k ) \mwC_{jk}\mwC_{j\ell}
       + \frac{3}{8\mwA^4N^3}\sum_{j,k=1}^N
                           (2\mwa_j+\mwa_k ) \mwB_j\mwC_{jk}
\displaybreak[0]\\
  \mwZ &= \exp\biggl(
       \frac{1}{4\mwA^2N^2}\sum_{j=1}^{N}\Im(\mwa_j)^2
     + \frac{15}{16\mwA^3N}\sum_{j=1}^{N}\Im(\mwB_j)^2 \\
     &\kern7mm{}
     + \frac{3}{8\mwA^3N^2}\sum_{j,k=1}^{N}\Im(\mwB_j)\Im(\mwC_{jk})
     + \frac{1}{16\mwA^3N^3}\sum_{j,k,\ell=1}^N \Im(\mwC_{jk})\Im(\mwC_{j\ell})
    \biggr).
\quad\qedsymbol
\end{align*}
\end{lemma}

\smallskip

\noindent We use the notation
$\R^c$ for the complement of a region $\R$.  To evaluate
the integral $I(m,n)$ defined in (\ref{Idef2}), we proceed as
follows:
$$
I(m,n) = \int_{\R}F + \int_{\R^c}F = \int_{\R'}F + O(1)\int_{\R^c}|F|,
~~~~ \R' \supseteq \R.
$$
The region $\R\subseteq [-\pi,+\pi]^{m+n}$
is defined below in (\ref{RDef}).  The larger containing
region $\R'$ is defined in (\ref{RprimeDef}).
The asymptotic value of the integral $\int_{\R'}F$ is obtained by
substituting~$\R'$ for the variable $\R''$ in equation (\ref{IvsJ}),
and applying the main result of this section,
Theorem~\ref{Jintegral}.
In the next section we show that for a certain constant $c_5>0$
$$
\int_{\R^c}|F| = O(e^{-c_5\min(m^{2\eps},n^{2\eps})/\log n})\int_{\R'}F.
$$
(Again, the quantity $\eps$ is a small positive constant arising in 
the definition (\ref{RDef}) of the region $\R$.)
This is the complete summary of how we shall evaluate
$I(m,n)$, and now we may proceed to the technical
details.

For any region $\R''$ we set
$$
I_{\R''}(m,n) = \int_{\R''}F.
$$
In the integral
(\ref{Idef2}), it is convenient sometimes to think of $\theta_j$,
$\phi_k$ as points on the unit circle.  We wish to define 
``averages'' of the angles $\theta_j$, $\phi_k$.  To do this cleanly
we make the following definitions, as in \cite{CM}.  Let $C$ be the
ring of real numbers modulo $2\pi$, which we can interpret as points
on a circle in the usual way.  Let $z$ be the canonical mapping
from $C$ to the real interval $(-\pi,\pi]$.  An {\it open half circle}
is $C_t=(t-\pi/2,t+\pi/2)\subseteq C$ for some~$t$.  Now define
$$
\Chat^N = \{\bix=(x_1,\dots,x_N) \in C^N | x_1,\dots,x_N \in C_t 
{\rm ~for~some~} t \in \Reals\}.
$$
If $\bix=(x_1,\dots,x_N) \in C_0^N$ then define
$$
\xav = z^{-1}\left( \frac{1}{N} \sum_{j=1}^N z(x_j)\right).
$$
More generally, if $\bix\in C_t^N$ then define $\xav =
t+\overline{(x_1-t,\dots,x_N-t)}$.  The function $\bix\mapsto\xav$
is well-defined and continuous for $\bix\in\Chat^N$.

For some sufficiently
small $\eps>0$, let $\R$ denote the set of vector pairs
$\bit,\bip\in \Chat^m\times\Chat^n$ such that
\begin{equation}\label{RDef}
\begin{split}
|\tb+\pb| &\le \tbd (mn)^{-1/2+2\eps} \\
{\rm region~} \R: ~~~~~~~~~~~~~~
|\th_j| &\le \tbd n^{-1/2+\eps}, 1 \le j \le m \\
|\ph_k| &\le \tbd m^{-1/2+\eps}, 1 \le k \le n,
\end{split}
\end{equation}
where $\th_j = \t_j-\tb$ and $\ph_k = \p_k-\pb$.
In this definition, values are considered in~$C$.

\medskip

Let $\bithat=(\th_1,\dots,\th_{m-1})$,
$\biphat=(\ph_1,\dots,\ph_{n-1})$, and define $T_1$ to
be the transformation
$T_1(\bithat,\biphat,\mu,\delta)=(\bit,\bip)$ given by
\[ \mu = \tav + \pav,
  \qquad  \delta = \tav -
\pav,\]
together with $\that_j=\t_j-\tav$
($1\leq j\leq m-1$) and $\phat_k=\p_k-\pav$ ($1\leq k\leq n-1$).
We also define the 1-many transformation $T_1^*$ by
\[
    T_1^*(\bithat,\biphat,\mu) =
    \bigcup_\delta \, T_1(\bithat,\biphat,\mu,\delta).
\]

After applying the transformation $T_1$ to $I_{\R}(m,n)$,
the new integrand is easily seen to be independent of
$\delta$, so we can multiply by the range of $\delta$ and remove
it as an independent variable.  Therefore, we can continue with
an $(m{+}n{-}1)$-dimensional integral over the region $\S$ defined by
$\R=T_1^*(\S)$.   More generally, if
$\S''\subseteq [-\tfrac12\pi,\tfrac12\pi]^{m+n-2}\times[-2\pi,2\pi]$
and $\R''=T_1^*(\S'')$, we have
\begin{equation}\label{IvsJ}
I_{\R''}(m,n) =  2\pi m n \int_{\S''} G(\bithat,\biphat,\mu)
     \, d\bithat  d\biphat d\mu,
\end{equation}
where $G(\bithat,\biphat,\mu)
=F\(T_1(\bithat,\biphat,\mu,0)\)$.
The factor $2\pi mn$ combines the range of $\delta$,
which is $4\pi$, and the Jacobian of $T_1$, which is~$mn/2$.

Note that  $\S$ is defined by virtually
the same inequalities as define~$\R$.  The first inequality
is now $\abs{\mu}\leq \tbd (mn)^{-1/2 + 2\eps}$ and the bounds on
\[ \that_m = - \sum_{j=1}^{m-1} \that_j~\text{ and }~
    \phat_n = - \sum_{k=1}^{n-1} \phat_k
\]
still apply even though these are no
longer variables of integration.

\medskip

Our main result in this section is the following.
\begin{thm}\label{Jintegral}
Define $A = \frac12\l(1+\l)$.
Under the conditions of Theorem~\ref{thm:main},  there is a
region $\S'\supseteq\S$ such that
\begin{align*}
 \int_{\S'} G(\bithat,\biphat,\mu)
     \, d\bithat  d\biphat  d\mu
     &= (mn)^{-1/2} \Bigl(\frac{\pi}{Amn }\Bigr)^{\!1/2}
                    \Bigl(\frac{\pi}{An  }\Bigr)^{(m-1)/2}
                    \Bigl(\frac{\pi}{Am  }\Bigr)^{(n-1)/2} \\
     &\kern-18mm {}\times \exp\biggl(
           \frac12
         - \frac{1+2A}{24A}\Bigl(\frac mn+\frac nm\Bigr)
         + O(n^{-b})  \biggr).
\end{align*}
\end{thm}

\begin{proof}
For $x$ real and $|\l x|$ small,
$$
(1-\l(e^{ix}-1))^{-1} = \exp\(
\l ix - Ax^2 - iA_3x^3 + A_4 x^4 +\Oh\((1+\l)^5|x|^5\) \)
$$
with
$$
A=\tfrac{1}{2}\l(1+\l),\;  A_3=\tfrac{1}{6}\l(1+\l)(1+2\l),\;
 A_4=\tfrac{1}{24}\l(1+\l)(1+6\l+6\l^2).
$$
We are about to state an estimate for $G$, and some further
estimates (\ref{EBounds}) in a moment, all of which hold uniformly
in the region $\S$.  Let us alert the reader that we shall be defining
a larger region, see (\ref{LargeR}) below, which contains $\S$ and
in which all of these estimates continue to hold true.
Uniformly in $\S$, where all
$(1+\l)|\mu+\that_j+\phat_k|$ are small,
\begin{align*}
G &= \exp\Bigl\{ -A\sumsum(\mu+\that_j+\phat_k)^2
                    - iA_3\sumsum(\mu+\that_j+\phat_k)^3
\\
    & \qquad\qquad{}
                   + A_4\sumsum(\mu+\that_j+\phat_k)^4
                   + \Oh\((1+\l)^5\,\sumsum|\mu+\that_j+\phat_k|^5\,\)
             \Bigr\}.
\end{align*}
Here and below, the undelimited summation over $j,k$ runs
over $1\le j\le m$, $1\le k\le n$, and we
continue to use the abbreviations
$\th_m=-\sum_{j=1}^{m-1}\th_j$ and
$\ph_n=-\sum_{k=1}^{n-1}\ph_k$. 

\smallskip

We now proceed to a second 
change of variables, $(\bithat,\biphat,\mu)=T_2(\bis,\bitau,\mu)$ given
by
$$
\that_j=\sigma_j+c\mu_1, ~~ \phat_k=\tau_k+d\nu_1,
$$
where, for $1\le h \le 4$, $\mu_h$ and $\nu_h$ denote the
power sums $\sum_{j=1}^{m-1}\sigma_j^h$ and
$\sum_{k=1}^{n-1}\tau_k^h$, respectively.
The scalars $c$ and $d$ are chosen to eliminate the
second-degree cross-terms $\sigma_{j_1}\sigma_{j_2}$
and $\tau_{k_1}\tau_{k_2}$, and thus diagonalize the
quadratic in $\bis=(\sigma_1,\dots,\sigma_{m-1})$
and $\bitau=(\tau_1,\dots,\tau_{n-1})$.  Suitable choices for $c,d$
are
$$
c = -\frac{1}{m+m^{1/2}}, ~~ d = -\frac{1}{n+n^{1/2}},
$$
and we find the following:
\begin{align}
\sumsum(\mu+\that_j+\phat_k)^2 &= 
     mn\mu^2
 + n\mu_2
 + m\nu_2
\notag
\\[0.5ex]
\sumsum(\mu+\that_j+\phat_k)^3 &= 
   3\mu(n\mu_2+m\nu_2)
 + n(\mu_3+3c\mu_2\mu_1)
 + m(\nu_3+3d\nu_2\nu_1)
\notag
\\
 & ~~~ + \tO((1+\l)^{-3}n^{-1/2})
\label{EBounds}
\\[0.5ex]
\sumsum(\mu+\that_j+\phat_k)^4 &=  
      6\mu_2\nu_2
    + n\mu_4
    + m\nu_4
    + \tO((1+\l)^{-4}n^{-1/2})
\notag
\\[0.5ex]
\sumsum|\mu+\that_j+\phat_k|^5
   &= \tO((1+\l)^{-5}n^{-1/2})\,.
\notag
\end{align} 
The Jacobian of the matrix $T_2$ is $(mn)^{-1/2}$,
and so
\begin{equation}
\int_{\S}G = (mn)^{-1/2}   \int_{T_2^{-1}(\S)} E_1,
\label{JEOne}
\end{equation}
where $E_1=\exp(L_2+\tO(n^{-1/2}))$, and
\begin{equation}\label{L2value}
\begin{split}
L_2 = 
-    A   n     \mu_2  -    A   m     \nu_2
-    A   m n \mu^2
-  i A_3 n     \mu_3  -  i A_3 m     \nu_3
- 3i A_3 n     \mu   \mu_2  - 3i A_3 m     \mu   \nu_2\\
- 3i A_3 c         n \mu_2  \mu_1
- 3i A_3 d         m \nu_2  \nu_1
+    A_4 n     \mu_4   +    A_4 m     \nu_4
+ 6  A_4 \mu_2 \nu_2.
\end{split}
\end{equation}

Define the regions $\Q$, $\M$, $\S'$, and $\R'$ by
\begin{align*}
\Q &= \bigl\{\,|\sigma_j| \le 2(1+\l)^{-1}n^{-1/2+\eps},
             \;j=1,\ldots,m{-}1\,\bigr\} \\
&\qquad\qquad{}\cap
      \bigl\{\,|\tau_k| \le 2(1+\l)^{-1}m^{-1/2+\eps},
             \;k=1,\ldots,n{-}1\,\bigr\},\\
&\qquad\qquad{}\cap
\bigl\{\,|\mu|\le 2(1+\l)^{-1}(mn)^{-1/2+2\eps}\,\bigr\}
\\[0.3ex]
\M &= \bigl\{\,|\mu_1| \le (1+\l)^{-1}m^{1/2}n^{-1/2+2\eps}\,\bigr\}
    \cap \bigl\{\,|\nu_1| \le (1+\l)^{-1}n^{1/2}m^{-1/2+2\eps}\,\bigr\},
\end{align*}
\begin{equation}\label{SprimeDef}
\S' = T_2(\Q\cap\M),
\end{equation}
and
\begin{equation}\label{RprimeDef}
\R' = T_1^{*}(\S').
\end{equation}
Summing for $1 \le j \le m-1$ the equation $\that_j=\sigma_j+c\mu_1$,
and inserting the value of $c$, we find
$$
m^{-1/2} \, \mu_1 = \sum_{j=1}^{m-1} \that_j.
$$
In the region $\S$ we have $|\sum_{j=1}^{m-1} \that_j| \le 
(1+\l)^{-1}n^{-1/2+\eps}$, and so in $T_2^{-1}(\S)$ we have
$$
|\mu_1| \le (1+\l)^{-1}m^{1/2}n^{-1/2+\eps}.
$$
Using this, and the dual inequality for $|\nu_1|$, the reader can check that
$$
T_2^{-1}(\S) \subseteq \Q \cap \M.
$$

It will be convenient if we can apply the bounds~\eqref{EBounds}
throughout the expanded region $\S'=T_2(\Q\cap\M)$, rather than
only in~$\S$.
To see that this is valid, note that the
calculations leading to these bounds hold equally well
if the coefficient of $\eps$ in an exponent of $m$ or $n$ is
increased; or, if an assumption 
such as $|\mu| \le \tbd (mn)^{1/2+2\eps}$ is made more permissive
by a multiplicative constant: $|\mu| \le 2 \tbd (mn)^{1/2+2\eps}$.
Therefore, it suffices to note that $S'$ lies inside the
region
\begin{align}
|\mu|&\le 2 \tbd (mn)^{-1/2+2\eps} \notag \\
|\th_j|    &\le 3 \tbd   n ^{-1/2+ \eps}, ~ 1 \le j \le m-1 \notag \\
|\th_m|    &\le   \tbd   n ^{-1/2+2\eps}, & \label{LargeR} \\
|\ph_k|    &\le 3 \tbd   m ^{-1/2+ \eps}, ~ 1 \le k \le n-1 \notag \\
|\ph_n|    &\le   \tbd   m ^{-1/2+2\eps}. \notag
\end{align}

Define $E_2=\exp(L_2)$.
We have shown that in the region $\Q\cap\M$ the integrand $E_1$
satisfies $E_1=E_2\(1+\tO(n^{-1/2})\)$.
We can
approximate our integral as follows:
\begin{align}
\int_{\Q \cap \M} E_1
  &= \int_{\Q \cap \M} E_2 + \tO(n^{-1/2})\int_{\Q \cap \M} \abs{E_2} \notag\\[0.4ex]
  &= \int_{\Q \cap \M} E_2 + \tO(n^{-1/2})\int_{\Q} \,\abs{E_2} \notag\\[0.4ex]
  &= \int_{\Q} E_2
        + O(1)\int_{\Q\cap \M^c} \abs{E_2} + \tO(n^{-1/2})\int_{\Q} \,\abs{E_2}.
                  \label{recipe}
\end{align}
It suffices to estimate each of the three integrals
in the last line of (\ref{recipe}).

We first compute the integral of $E_2$ over $\Q$.
We proceed in three stages, starting with integration
with respect to $\mu$.
For the latter, we can use the formula
\[
 \int\limits_{-(mn)^{-1/2+2\eps}}^{(mn)^{-1/2+2\eps}}
   \kern-2mm\exp\(-Amn\mu^2 - i\beta\mu\) d\mu
 = \Bigl(\frac{\pi}{Amn}\Bigr)^{1/2}
   \exp\biggl( -\frac{\beta^2}{4Amn} +O(n^{-1}) \biggr),
\]
provided $\beta=o(A(mn)^{1/2+2\eps})$.  In our case,
$\beta=3A_3(n\mu_2+m\nu_2)$, which is small enough
because $m=O(n\log n)$
and $n=O(m\log m)$.
Integration over $\mu$ contributes
\begin{equation}\label{muint}
\Bigl(\frac{\pi}{Amn}\Bigr)^{1/2}
 \exp\biggl(\frac{-9A_3^2(n\mu_2+m\nu_2)^2}{4Amn}
+O(n^{-1}) \biggr).
\end{equation}

The second step is to
integrate with respect to $\bis$ the integrand
\begin{align}\label{sigbits}
\begin{split}
   \exp\biggl(\! &- An\mu_2
         - \frac{9A_3^2n}{4Am}\mu_2^2
          - i A_3 n\mu_3
               - 3i A_3 cn \mu_1\mu_2  \\
       &~{}+ \Bigl(6A_4-\frac{9A_3^2}{2A}\Bigr)\mu_2\nu_2 + A_4n\mu_4
       + O(n^{-1})\biggr).    
\end{split}
\end{align}
This is accomplished by an appeal to Lemma~\ref{MW3}.
We must take the $N$ of that lemma equal to $m-1$, the
number of variables.  This dictates that the limits of
integration be $\pm (m-1)^{1/2+\hat{\eps}}$, but our
limits, based on the definition of $\Q$, are
$\pm 2(1+\l)^{-1}n^{-1/2+\eps}$.  Thus, we make a
change of scale,
$$
\sigma_j \leftarrow \frac{\sigma_j}{(1+\lambda)(n/m)^{1/2}}
$$
before integrating.  This change of scale will introduce
a multiplicative factor $(1+\lambda)^{-(m-1)}(m/n)^{(m-1)/2}$ into our
evaluation of the integral.
In the terminology of that lemma, we have
$N = m-1$, $\delta(N)=O(n^{-1})$, 
$\eps'=\tfrac32\eps$, $\eps''=\tfrac53\eps$, $\eps'''=3\eps$,
$\bar\eps=6\eps$, and $\hat\eps(N)=\eps+o(1)$ is defined
by $2m^{-1/2+\eps}=(m-1)^{-1/2+\hat\eps}$.
Furthermore, taking account of the scale, we find
\begin{align*}
\mwA &= \frac{A}{(1+\l)^2}\frac{m}{m-1}, &
\mwa_j &=  \frac{12AA_4-9A_3^2}{2A(1+\l)^2} \frac{m}{n}\nu_2, &&\\
\mwB_j &= -\frac{iA_3}{(1+\l)^3}\frac{m^{1/2}}{n^{1/2}}\frac{m}{m-1},  &
\mwC_{jj'} &= -\frac{3iA_3}{(1+\l)^3} cn \frac{m^{3/2}}{n^{3/2}}, && \\
\mwE_j &= \frac{A_4}{(1+\l)^4} \frac{m}{n} \frac{m}{m-1}, &
\mwF_{jj'} &= -\frac{9A_3^2}{4A(1+\l)^4} \frac{m}{n}, & \mwJ_j &= 0.
\end{align*}
We can take $\Delta=\tfrac45$, and
calculate that
\begin{gather}
 \frac{3}{4\mwA^2N}\sum_{j=1}^N \mwE_j
    + \frac{1}{4\mwA^2N^2}\sum_{j,j'=1}^N \mwF_{jj'} = \frac{m}{n}
  \biggl(\frac{3A_4}{4A^2}-\frac{9A_3^2}{16A^3}\biggr) + \tO(n^{-1}) \notag\\
 \frac{15}{16\mwA^3N}\sum_{j=1}^{N}\mwB_j^2
     + \frac{3}{8\mwA^3N^2}\sum_{j,j'=1}^{N}\mwB_j\mwC_{jj'}
  + \frac{1}{16\mwA^3N^3}\sum_{j,j'\!,j''=1}^N \mwC_{jj'}\mwC_{jj''}
          =  -\frac{3A_3^2m}{8A^3n} + \tO(n^{-1}) \notag\\
          \frac{1}{2\mwA N}\sum_{j=1}^{N}\mwa_j
     + \frac{1}{4\mwA^2N^2}\sum_{j=1}^{N}\mwa_j^2
      = 
    \frac{m}{n}\biggl(\frac{3A_4}{A}-\frac{9A_3^2}{4A^2}\biggr)\nu_2
              +\tO(n^{-1})  
     \label{taubit}\\
    \mwZ = Z_1 = \exp\biggl( \frac{3A_3^2m}{8A^3n} + \tO(n^{-1}) \biggr)
         = O(1) \exp\biggl(\frac {(1+2\lambda)^2 m}{24 An}\biggr).\notag
\end{gather}

After checking that  $\Theta_2 = \tO(n^{-1})=o(m^{-4/5}Z_1)$,
we conclude that integration with respect to
$\bis$ contributes a $\bitau$-free factor
\begin{equation}\label{contrib1}
 \Bigl(\frac{\pi}{An}\Bigr)^{\!(m-1)/2}
 \exp\biggl(
  \Bigl(\frac{3A_4}{4A^2}-\frac{15A_3^2}{16A^3}\Bigr) \frac{m}{n}
    + O(m^{-4/5}Z_1)
    \biggr).
\end{equation}
By the conditions of Theorem~\ref{thm:main}, $Z_1\le n^{1/5}$,
so $m^{-4/5}Z_1=o(1)$ as required by Lemma~\ref{MW3}.

Finally, we need to integrate over $\bitau$.  Collecting the
remaining terms from \eqref{L2value}, and the terms
involving $\bitau$ from
(\ref{muint}) and (\ref{taubit}), we have an
integrand equal to
\begin{align*}
&
\exp\biggl(
-Am\nu_2
+\Bigl(\frac{3A_4m}{An}-\frac{9A_3^2m}{4A^2n}\Bigr)\nu_2
-\frac{9A_3^2m}{4An}\nu_2^2
\\
& \kern11mm{}
+A_4m\nu_4  -iA_3m\nu_3 -3iA_3dm\nu_2\nu_1
  + O(m^{-4/5}Z_1)
\biggr).
\end{align*}

Again we use Lemma~\ref{MW3}.  This time, the factor due to
scaling is $(1+\l)^{-(n-1)}(n/m)^{(n-1)/2}$, and 
$N = n-1$, $\delta(N)=O(m^{-4/5}Z_1)$,
$\eps'=\tfrac32\eps$, $\eps''=\tfrac53\eps$, $\eps'''=3\eps$,
$\bar\eps=4\eps$, and $\hat\eps(N)=\eps+o(1)$, as defined by
$2n^{-1/2+\eps}=(n-1)^{-1/2+\hat\eps}$.
The substitution table this time reads
\begin{align*}
\mwA &= \frac{A}{(1+\l)^2}\frac{n}{n-1}, &
\mwa_k &=  \frac{12AA_4-9A_3^2}{4A^2(1+\l)^2},  &&\\
\mwB_k &= -\frac{iA_3}{(1+\l)^3}\frac{n^{1/2}}{m^{1/2}}\frac{n}{n-1},  &
\mwC_{kk'} &= -\frac{3iA_3}{(1+\l)^3} dm \frac{n^{3/2}}{m^{3/2}}, && \\
\mwE_k &= \frac{A_4}{(1+\l)^4} \frac{n}{m} \frac{n}{n-1}, &
\mwF_{kk'} &= -\frac{9A_3^2}{4A(1+\l)^4} \frac{n}{m}, & \mwJ_k &= 0.
\end{align*}
This time we take $\Delta_2>4/5$ so that
$n^{-\Delta_2}=o(\delta(N))$; calculations similar to the
previous case lead to
\begin{gather*}
 \frac{3}{4\mwA^2N}\sum_{k=1}^N \mwE_k
    + \frac{1}{4\mwA^2N^2}\sum_{k,k'=1}^N \mwF_{kk'} = \frac{n}{m}
  \biggl(\frac{3A_4}{4A^2}-\frac{9A_3^2}{16A^3}\biggr) + \tO(n^{-1}) 
  \\
 \frac{15}{16\mwA^3N}\sum_{k=1}^{N}\mwB_k^2
     + \frac{3}{8\mwA^3N^2}\sum_{k,k'=1}^{N}\mwB_j\mwC_{kk'}
   + \frac{1}{16\mwA^3N^3}\sum_{k,k'\!,k''=1}^N \mwC_{kk'}\mwC_{kk''}\\
        {\kern 8cm} =  -\frac{3A_3^2n}{8A^3m} + \tO(n^{-1}) \\
          \frac{1}{2\mwA N}\sum_{k=1}^{N}\mwa_k
     + \frac{1}{4\mwA^2N^2}\sum_{k=1}^{N}\mwa_k^2
      =
          - \frac{9A_3^2}{8A^3} + \frac{3A_4}{2A^2} + \tO(n^{-1}) \\
    \mwZ = Z_2 = \exp\biggl( \frac{3A_3^2n}{8A^3m} + \tO(n^{-1}) \biggr)
         = O(1) \exp\biggl(\frac {(1+2\lambda)^2 n}{24 Am}\biggr).
\end{gather*}

Again $\Theta_2=\tO(n^{-1})$, implying that $\Theta_2=o(m^{-4/5}Z_1Z_2)$.
Including the contributions
from \eqref{muint} and \eqref{contrib1}, we obtain
\begin{equation}\label{QE2a}
\begin{split}
 \int_\Q E_2
  &= \Bigl(\frac{\pi}{Amn}\Bigr)^{1/2}
        \Bigl(\frac{\pi}{An}\Bigr)^{\!(m-1)/2}
        \Bigl(\frac{\pi}{Am}\Bigr)^{\!(n-1)/2} \\
  &{\kern-5mm}\times\exp\biggl(  - \frac{9A_3^2}{8A^3} + \frac{3A_4}{2A^2}
         +\Bigl(\frac mn+\frac nm\Bigr)
         \Bigl( \frac{3A_4}{4A^2}-\frac{15A_3^2}{16A^3}\Bigr) + O(m^{-4/5}Z_1Z_2) \biggr).
\end{split}
\end{equation}

Since
$$
Z_1Z_2=O(1)  \exp\left( \frac{3A_3^2}{8A^3} \left(
                                                    \frac{m}{n} + \frac{n}{m}
                                           \right)
                \right), 
$$
and since 
$$
\frac{3A_3^2}{8A^3} = \frac{(1+2\l)^2}{24A},
$$
it follows from the main hypothesis (\ref{Hyp})
of Theorem~\ref{thm:main} that $Z_1Z_2 = O(n^{6a/5})$.
The condition $a+b<1/2$ implies $-4/5+6a/5<-b-1/5$; hence,
substituting the values of $A, A_3, A_4$, we conclude that
\begin{equation}\label{QE2}
\begin{split}
 \int_\Q E_2
  &= \Bigl(\frac{\pi}{Amn}\Bigr)^{1/2}
        \Bigl(\frac{\pi}{An}\Bigr)^{\!(m-1)/2}
        \Bigl(\frac{\pi}{Am}\Bigr)^{\!(n-1)/2} \\
  &\qquad{}\times\exp\biggl(  \frac12
         - \frac{1+2A}{24A}\Bigl(\frac mn+\frac nm\Bigr)
      + O(n^{-b-1/5}) \biggr).
\end{split}
\end{equation}

\medskip
We next infer a estimate of $\int_\Q\, \abs{E_2}$.  The calculation that
lead to \eqref{QE2a} remains valid if we set
$A_3$
to zero, which is the same as replacing $L_2$ by its real part.
Since $\abs{E_2} = \exp\(\Re(L_2)\)$, this gives
\begin{align}
  \int_\Q\, \abs{E_2}
    &= \exp\biggl(\frac{(1+2\lambda)^2}{8A}\Bigl( 1 +
         \frac{5n}{6m} + \frac{5m}{6n}\Bigr)+ o(1)\biggr)\int_Q E_2 \notag\\
     &= O(n^a)\int_Q E_2 \label{absE2}
\end{align}
under the assumptions of Theorem~\ref{thm:main}.
The third term of (\ref{recipe}) can now be identified:
\begin{equation}\label{OEbit}
   \tO(n^{-1/2})  \int_\Q\, \abs{E_2} = \tO(n^{-1/2+a}) \int_Q E_2
        = O(n^{-b})  \int_Q E_2\,.
\end{equation}

Finally, we consider the second term of \eqref{recipe}, namely
\[\int_{\Q\cap \M^c} \abs{E_2},\]
which we will bound as a fraction of $\int_\Q\, \abs{E_2}$ using
a statistical technique.  The following is a well-known result of
Hoeffding~\cite{hoeffding}.

\begin{lemma}\label{hoeffding}
  Let $X_1,X_2,\ldots,X_N$ be independent random
variables such that $EX_i=0$ and $\abs{X_i}\le M$ for all~$i$.
Then, for any $t\ge 0$,
\[\Prob\biggl(\,\sum_{i=1}^N X_i \ge t\biggr)
    \le \exp\biggl(-\frac{t^2}{2NM^2}\biggr).\]
\end{lemma}

Now consider $\abs{E_2} = \exp\(\Re(L_2)\)$.
Write $\M=\M_1\cap\M_2$, where
$\M_1=\{\, \abs{\mu_1}\le m^{1/2} n^{-1/2+2\eps}\,\}$
and $\M_2=\{\, \abs{\nu_1}\le n^{1/2} m^{-1/2+2\eps}\}$.
For fixed values
of $\mu$ and $\bis$, $\Re(L_2)$ separates over
$\tau_1,\tau_2,\ldots,\tau_{n-1}$ and therefore, apart from
normalization,  it is the joint density of independent
random variables $X_1,X_2,\ldots,X_{n-1}$ which satisfy
$EX_k=0$ (by symmetry) and $\abs{X_k}\le 2(1+\l)^{-1}m^{-1/2+\eps}$
(by the definition of $\Q$).  By Lemma~\ref{hoeffding}, the
fraction of the integral over $\bitau$ (for fixed $\mu,\bis$)
that has
$\nu_1\ge (1+\l)^{-1}n^{1/2}m^{-1/2+2\eps}$ is at most
$\exp(-m^{2\eps}/8)$.  By symmetry, the same bound holds
for $\nu_1\le -(1+\l)^{-1}n^{1/2}m^{-1/2+2\eps}$.  Since these bounds
are independent of $\mu$ and $\bis$, we have
\[\int_{\Q\cap\M_2^c} \,\abs{E_2}
   \le 2\exp(-m^{2\eps}/8)\int_\Q\,\abs{E_2}.\]
By the same argument,
\[\int_{\Q\cap\M_1^c} \,\abs{E_2}
   \le 2\exp(-n^{2\eps}/8)\int_\Q\,\abs{E_2}.\]
Therefore we have in total that
\begin{align}\label{corners}
\int_{\Q\cap\M^c} \,\abs{E_2}
   &\le 2\(\exp(-m^{2\eps}/8)+\exp(-n^{2\eps}/8)\) \int_\Q\,\abs{E_2}
      \notag \\
   &\le  O(n^{-b})  \int_Q E_2,
\end{align}
using again (\ref{QE2}).  Applying \eqref{recipe} with \eqref{QE2},
\eqref{OEbit} and~\eqref{corners}, we find that $\int_{\Q\cap\M} E_1$
is given by~\eqref{QE2} with the error term replaced by~$O(n^{-b})$.
Multiplying by the Jacobian of the transformation~$T_2$,
we find that Theorem~\ref{Jintegral} is proved for $\S'$
given by~(\ref{SprimeDef}).
\end{proof}

\nicebreak
\section{Concentration of the integral}\label{section:Boxing}
Recall that $F(\bit,\bip)$ is the integrand in equation (\ref{Idef2})
defining $I(m,n)$, and that $\R$ is the region given by (\ref{RDef}).
In the previous section we estimated the integral of $F(\bit,\bip)$
over a particular superset $\R'\supseteq\R$.
In this section we show that the integral of $F(\bit,\bip)$
outside $\R$ is negligible in comparison if $\lambda$ is
polynomially bounded.
Larger values of $\lambda$ will be handled in the
following section.

\begin{thm}\label{boxing}
Suppose that $m,n \rightarrow \infty$
in such a way that~\eqref{Hyp} holds and $\lambda=n^{O(1)}$.
Define $I_0$ by
\begin{equation}\label{ThmAssumption}
 I_0 = (mn)^{1/2}
                    \Bigl(\frac{\pi}{Amn }\Bigr)^{ 1/2}
                    \Bigl(\frac{\pi}{An  }\Bigr)^{ (m-1)/2}
                    \Bigl(\frac{\pi}{Am  }\Bigr)^{ (n-1)/2}
              \! \exp\biggl( 
- \frac{1+2A}{24A}\Bigl(\frac mn+ \frac nm\Bigr)
        \biggr).
\end{equation}
Then, for sufficiently small $\eps>0$,
$$
\int_{\R^c} |F| = O(e^{-n^{\eps}})I_0.$$
\end{thm}

%
We begin with two technical lemmas whose proofs are omitted.

\begin{lemma}\label{fbnd}
The absolute value of the integrand of $I(m,n)$ is
$$|F(\bit,\bip)| = \prod_{j,k} f(\t_j + \p_k),$$
where
$$f(z) = \(1+4A(1-\cos z))^{-1/2}.$$
Moreover, for all real $z$ with $\abs z\le \tfrac1{10}(1+\l)^{-1}$,
$$0\le f(z) \le \exp\( -A z^2 + (\tfrac1{12} A+A^2) z^4\).
    \quad\qedsymbol$$
\end{lemma}

\begin{lemma}\label{ibnd}
Define $N = \lceil 6000 (1+\l)\rceil$, $\d=2\pi/N$
and $g(x)=-Ax^2+(\frac94A+27A^2)x^4$.
Then, uniformly for $\l>0$, $K\ge 1$,
$$\int_{-30\d}^{30\d} \exp\( K g(x) \)\,dx
 \le \sqrt{\pi/(AK)}\, \exp\( O(K^{-1}+(AK)^{-1})\).
    \quad\qedsymbol$$
 \end{lemma}


\begin{proof}[Proof of Theorem \ref{boxing}]
Let $N$
and $\delta$ be as given in Lemma 4. Define the
region $\A$ to be the set of those $(\bit,\bip)$ such that
$$
\cos(\theta_j+\phi_k) \le \cos \delta
$$
for at least $\tfrac13\min(mn^{\eps},m^{\eps}n)$ pairs
$(j,k)$. Define $x_0,x_1,\ldots,x_{N-1}$ by $x_\ell=\ell\delta$.

If $X\subseteq (-\pi,\pi]$, we denote by $\Ntheta(X)$ the number
of values of $j$ such that $\theta_j\in X$, and similarly
define $\Nphi(X)$.
Define region $\R_1(\ell)$ to be the set of
those $(\bit,\bip)$ such that
$\Ntheta([x_\ell-4\delta,x_\ell+4\delta])\ge m-m^\eps$ and
$\Nphi([-x_\ell-4\delta,-x_\ell+4\delta])\ge n-n^\eps$.

Let $\U=\bigcup_{\ell=0}^{N-1}\R_1(\ell)$.  The proof of the
theorem consists in proving these three relations:
\begin{align}
\A\cup\U &= [-\pi,\pi]^{m+n}\label{R1}\\[0.8ex]
\int_{\A}|F| &= O(e^{-n})I_0\label{R2}\\[0.5ex]
\int_{\U\cap\R^c}|F| 
 &= O(e^{-n^\eps})I_0.\label{R3}
\end{align}

To begin the proof of~\eqref{R1}, we show that
any point $(\bit,\bip)$ for which
$\Ntheta([x_\ell-\delta,x_\ell+\delta]\ge m^\eps$
belongs to $\A\cup\U$.
Indeed, if such a point does
not belong to $\A$, then it must have
$\Nphi([-x_\ell-2\delta,-x_\ell+2\delta])\ge \tfrac23 n$.
This in turn forces 
$\Ntheta([x_\ell-3\delta,x_\ell+3\delta])\ge m-m^\eps$,
which forces
$\Nphi([-x_\ell-4\delta,-x_\ell+4\delta])\ge n-n^\eps$.
In particular, $(\bit,\bip)\in\R_1(\ell)$.

To complete the proof of~\eqref{R1}, we show that
$(\bit,\bip)$ belongs to $\A$ if
$\Ntheta([x_\ell-\delta,x_\ell+\delta])\le m^\eps$
for all~$\ell$.
Let $a$ be a minimum-length interval $[x_\ell,x_{\ell'}]$
such that $\Ntheta(a)\ge m^\eps$.
Then $\Ntheta(a)\le 2m^\eps$ and the complementary interval
$\overline{a}$ has $\Ntheta(\overline a)\ge m-2m^\eps$.
If $b=[x_{\ell'+2},x_{\ell-2}]$ (a subinterval of $\overline a$),
then $\Ntheta(b)\ge m-4m^\eps$.
Thus, there are at least $m^\eps$
disjoint pairs $(j,j')$ with $\theta_j\in a$ and $\theta_{j'}\in b$.
(By disjoint we mean there are $2m^\eps$ distinct indices $j$
involved.)
Because of the $2\delta$ spaces between arcs $a$ and $b$,
for each $k$ and each pair $(j,j')$ at least one of
$\cos(\theta_j+\phi_k)$ or $\cos(\theta_{j'}+\phi_k)$ is
bounded above by $\cos\delta$.  This implies
that $(\bit,\bip)\in\A$, as claimed, and completes the proof
of~\eqref{R1}.

%

\smallbreak

We turn next to~\eqref{R2}. 
Since $A\delta^2=\Theta\(\l(1+\l)^{-1}\)$, Lemma 3 implies
\[
\abs{F(\bit,\bip)} \le \exp\(-c_1\l(1+\l)^{-1}\min(m^{\eps}n,mn^{\eps})\)
\]
for all $(\bit,\bip)\in\A$ and some $c_1>0$.
The volume of $\A$ is no more than $(2\pi)^{m+n}$, and 
$\l(1+\l)^{-1}>(\log n)^{-1}$ by~\eqref{Hyp}, so
\[
 \int_{\A}|F| 
   \le (2\pi)^{m+n}\exp\(-c_1(\log n)^{-1}\min(m^{\eps}n,mn^{\eps})\).
\]
{}From~\eqref{Hyp}, which implies that $A=O(\log n)$,
and the assumption that $\lambda=n^{O(1)}$, we have
that $I_0=\exp\(O(m\log n+n\log m)\)$.
Relation~\eqref{R2} follows.

\smallbreak

Finally we come to~\eqref{R3}.
For $(\bit,\bip)\in\R_1(\ell)$ we define
$S_0=S_0(\bit)$, $S_1=S_1(\bit)$ and $S_2=S_2(\bit)$
to be the set of indices $j$ such that $|\t_j-x_{\ell}|$
is: less than or equal to $4\d$,
in the interval $(4\d,5\d]$, and
larger than $5\d$, respectively.  The index sets 
$T_0=T_0(\bip)$, $T_1=T_1(\bip)$ and $T_2=T_2(\bip)$
are defined similarly.  Define $\R_1(\ell;m_2,n_2)$ to be
that subregion of $\R_1(\ell)$ for which
$|S_2|=m_2$ and $|T_2|=n_2$.
Define $\U(m_2,n_2)$ by
$$
\U(m_2,n_2) = \bigcup_{\ell=0}^{N-1}\R_1(\ell;m_2,n_2).
$$
We note that $m_2$ and $n_2$ vary
over the ranges $[0,m^{\eps}]$ and
$[0,n^{\eps}]$, respectively.
Define $\U_0=\U(0,0)$ and $\U_*=\U\setminus\U_0$.

Suppose $(\bit,\bip)\in\R_1(0;m_2,n_2)$.
If $|\t_j|$ and $|\p_k|$ are both less than or equal
to $5\d$, then Lemma 3 is applicable to $f(\t_j+\p_k)$.
If one of the two is less than or equal to $4\d$ and the
other exceeds $5\d$, then $\cos(\t_j+\p_k)\le\cos\d$.
Thus, there exists $c_2>0$ such that for
$(\bit,\bip)\in\R_1(0;m_2,n_2)$
\begin{equation}\label{R4}
f(\t_j+\p_k) \le 
  \begin{cases}
     \,\exp\(-A(\t_j+\p_k)^2+(\tfrac1{12}A+A^2)(\t_j+\p_k)^4\), \kern -32mm& \\
 	& \textrm{if }(j,k) \in (S_0\cup S_1)\times(T_0\cup T_1),\\[0.8ex]
           \,\exp\( -2c_2\l(1{+}\l)^{-1} \),
         & \textrm{if }(j,k) \in (S_0\times T_2)\cup (S_2\times T_0),\\[0.3ex]
     \,1,  & \textrm{otherwise}.
 \end{cases}
\end{equation}
Since $(\bit,\bip)\in\U$, we have
$|S_0|\ge m-m^{\eps}$ and $|T_0|\ge n-n^{\eps}$; thus,
the size of $(S_0\times T_2)\cup (S_2\times T_0)$ exceeds 
$\tfrac12(mn_2+nm_2)$.  Integrating the
upper bound on $|F|$ implied by (\ref{R4}), we find
\begin{equation}\label{R22}
\int_{\R_1(0;m_2,n_2)}\negthickspace\negthickspace \abs{F}
    \le (2\pi)^{m_2+n_2}
     \binom{m}{m_2}\binom{n}{n_2}
      \exp\(-c_2\l(1{+}\l)^{-1}(mn_2+m_2n)\)
       \, I'_2(m_2,n_2)
\end{equation}
with
\[
 I'_2(m_2,n_2) = \int_{-5\d}^{5\d}\!\!\!\cdots \int_{-5\d}^{5\d}
    \exp\Bigl(- A\sumpp(\t_j+\p_k)^2 +
		(\tfrac1{12}A+A^2)\sumpp(\t_j+\p_k)^4\Bigr)
      \, d\bit'' d\bip'',\notag
\]
where the double-primes denote restriction to $j\in S_0\cup S_1$
and $k\in T_0\cup T_1$.  
The factor $(2\pi)^{m_2+n_2}$ comes from
integrating over $\t_j$ for $j\in S_2$ and $\p_k$ for $k\in T_2$,
while the binomial coefficients account for the choices of
$S_2$ and $T_2$.

Set $m'=m-m_2,n'=n-n_2$. 
As implied by the notation, the value of $I'_2(m_2,n_2)$
is independent of which specific variables constitute the sets
$S_0\cup S_1$ and $T_0\cup T_1$.

Summing on $m_2+n_2\ge 1$
yields an upper bound for the integral of $|F|$ over the region
$\bigcup_{m_2+n_2\ge 1}\R_1(0;m_2,n_2)$.
Notice, however, that the transformation $\t_j\mapsto\t_j-x_\ell$,
$\p_k\mapsto\p_k+x_\ell$ is a volume-preserving bijection of
$\R_1(\ell;m_2,n_2)$ with $\R_1(0;m_2,n_2)$ which leaves
$F(\bit,\bip)$ invariant.  Thus, introducing an additional
factor of~$N$ to account for all values of~$\ell$, we have
\begin{align}
\int_{\U_*}|F|
  &  \,\le\,  N\negthickspace \sum_{m_2+n_2\ge 1}
\binom{m}{m_2}\binom{n}{n_2}(2\pi)^{m_2+n_2}
	&\notag \\
  &  \qquad{}\times
\exp\( -c_2\l(1+\l)^{-1}(mn_2+nm_2) \)
              \, I'_2(m_2,n_2). \label{R6}
\end{align}

We now analyze $I_2'(m_2,n_2)$ more closely.
Define
$\tav'=(m')^{-1}\sumpp\t_j$,
$\breve\t_j=\t_j-\tav'$ for~$j\in S_0\cup S_1$,
$\pav'=(n')^{-1}\sumpp\p_k$,
$\breve\p_k=\p_k-\pav'$ for~$k\in T_0\cup T_1$,
$\mu'=\pav'+\tav'$ and $\nu'=\pav'-\tav'$.
In terms of $\mu',\nu', \breve\t_j, \breve\p_k$
we have
$$ \sumpp(\t_j+\p_k)^2 =
   m'n'{\mu'}^2 + n'\sumpp_j\breve\t_j^2 
   + m'\sumpp_k \breve\p_k^2
$$
and
$$ \sumpp(\t_j+\p_k)^4 \le
   27m'n'{\mu'}^4 + 27n'\sumpp_j\breve\t_j^4 + 27m'\sumpp_k \breve\p_k^4.
$$
The latter follows from the inequality
$(x+y+z)^4 \le 27(x^4+y^4+z^4)$ valid for all $x,y,z$.
It follows that
\begin{equation}
\label{R5}
I'_2(m_2,n_2) \le 
\int_{-5\d}^{5\d}\!\!\!\cdots \int_{-5\d}^{5\d}
    \exp\Bigl(
           m'n' g(\mu') 
             + n'\sumpp_j g(\breve\t_j) + m'\sumpp_k g(\breve\p_k)
         \Bigr)  \, d\bit'' d\bip''.
\end{equation}
For the moment we are thinking of $\mu'$, $\nu'$, $\breve\t_j$,
and $\breve\p_k$ ($1\le j\le m', 1\le k\le n'$) as functions of the 
variables of integration.  Arbitrarily choose one element
$j_*\in S_0\cup S_1$ and one element $k_*\in T_0\cup T_1$.  Define
$S_3=(S_0\cup S_1)-\{j_*\}$ and $T_3=(T_0\cup T_1)-\{k_*\}$.  We claim
that the inequality (\ref{R5}) remains valid when
$\sumpp_j \equiv \sum_{j\in S_0\cup S_1}$ is replaced by 
$\sumpp_{j\neq j_*} \equiv \sum_{j\in S_3}$, and similarly for $k_*$.
To verify this assertion, note that since $|\t_j|\le 5\d$ for all
$j\in S_0\cup S_1$ it must be the case that $|\tav'|\le 5\d$, too, and so
$|\breve\t_{j_*}|\le 10\d$.  The claim now follows because, as is
readily checked, $g(x)\le 0$ for $|x|\le 10\d$.

As in the previous section, the next step is to make a
change of variables to $\mu'$, $\nu'$, $\breve\t_j$ ($j\in S_3$), 
and $\breve\p_k$ ($k\in T_3$). 
The region corresponding under this change of variables to
$[-5\d,5\d]^{m'+n'}$ is not exactly a product, but it is contained in the
product $[-10\d,10\d]^{m'+n'}$.  (The argument given a few lines earlier
to show $|\breve\t_{j_*}|\le 10\d$ applies equally well to
$|\breve\t_j|$, $j\in S_3$.)  Hence,
\begin{align*}
I'_2(m_2,n_2) \le O(\d)m'n' 
    \int_{-10\d}^{10\d}\!\!\cdots \int_{-10\d}^{10\d}
     \exp\Bigl(&  m'n' g(\mu') + n' \sum\nolimits_{j\in S_3} g(\breve\t_j)  \\
    &{}+ m' \sum\nolimits_{k\in T_3} g(\breve\p_k) \Bigr)\,
    d\breve\t_{j\in S_3}d\breve\p_{k\in T_3}\, d\mu'.
\end{align*}
Here, the factor $O(\d)m'n'$ comes from the integration of
$\nu'$ over $[-10\d,10\d]$ and the Jacobian $m'n'/2$ of the transformation.
Now the integral splits into a product of $m'+n'-1$
one-dimensional integrals.  Each of these factors can be
bounded by  Lemma~\ref{ibnd}, and we find that
\begin{align*}
  I'_2(m_2,n_2) &\le O(\d) m' n'
        \Bigl(\frac{\pi}{Am'n'}\Bigr)^{1/2}
   \,   \Bigl(\frac{\pi}{Am'}\Bigr)^{(n'-1)/2}
   \,   \Bigl(\frac{\pi}{An'}\Bigr)^{(m'-1)/2} \\
   &\qquad{}\times\exp\( o(1) + O(1+A^{-1})\(m'/n' + n'/m'\) \).
\end{align*}
By $m'\sim m,n'\sim n$ and (\ref{Hyp}), the $\exp(\cdots\)$ term
above equals $n^{O(1)}$.  Also,
$$
\Bigl(
         \frac{\pi}{Am'}
\Bigr)^{(n'-1)/2}
 = 
\Bigl(
         \frac{\pi}{Am}
\Bigr)^{(n-1)/2}
n^{O(m_2+n_2)},
$$
and similarly for the other terms.
Noting that $\d N=O(1)$, we find from~\eqref{R6} that
\[
\int_{\U_*}|F| \le
 \sum_{m_2+n_2\ge 1} n^{O(m_2+n_2)}
     \exp\( -c_3\l(1+\l)^{-1}(mn_2+nm_2)\) I_0.
\]
Since $\l=\Omega\((\log n)^{-1}\)$, the sum is dominated
by the terms with $m_2+n_2=1$.  We conclude that
\begin{equation}\label{R7}
\int_{\U_*}|F| = O(e^{-n^{1-\eps}}) I_0.
\end{equation}


\medskip

It remains to consider $\U_0\cap\R^c$.   As was argued before in obtaining
(\ref{R22}) we have in the $m_2=n_2=0$ case
$$
\int_{\R_1(0;0,0)\cap\R^c}|F|
 \le
\int_{[-5\d,5\d]^{m+n}\cap\R^c} \negthickspace
    \exp\Bigl(- A\sum(\t_j+\p_k)^2 +
                (\tfrac1{12}A+A^2)\sum(\t_j+\p_k)^4\Bigr)
      \, d\bit d\bip,
$$
where now the sums are over all $j,k$.  Observe that
the transformation
$\t_j\mapsto\t_j-x_{\ell},\p_k\mapsto\p_k+x_{\ell}$
takes $\R_1(\ell;0,0)\cap\R^c$ bijectively to $\R_1(0;0,0)\cap\R^c$,
since $\R,\R^c$ are invariant under this mapping.
Hence, the integral over $\U_0\cap\R^c$ of $|F|$ is no larger
than $N$ times the integral on the right side of the last displayed
relation.

The region of integration, $[-5\d,5\d]^{m+n}\cap\R^c$,
is contained in the union $\bigcup_{h=1}^{m+n+1}\P_h$,
where $\P_h$ equals the product $[-5\d,5\d]^{m+n}$ intersected with
that part of $\R^c$ where the $h$th inequality of definition
(\ref{RDef}) fails.  As was done in the $\U_*$ case,
we may bound $\int_{\P_h}\cdots$ by a product of $m+n-1$ 
one-dimensional integrals.  However, in the present situation,
when we make the change of variables to $(\mu',\nu',\breve\t_j,\breve\p_k)$,
the transformed region of integration, albeit still contained
in the product $[-10\d,10\d]^{m+n-1}$ (we are omitting $\nu'$ from
the discussion), has the additional property that one of the
latter intervals, the one corresponding to the index $h$,
is missing a neighborhood of size $(1+\l)^{-1}(mn)^{-1/2+2\eps}$,
$(1+\l)^{-1}n^{-1/2+\eps}$, or $(1+\l)^{-1}m^{-1/2+\eps}$ about 0.  
 
Throughout $[-10\d,10\d]$ we have $g(x)\le-Ax^2/2$.  Because
of the inequality  
$$
\int_{(1+\l)^{-1}K^{-1/2+\alpha}}^{\infty}e^{-KAx^2/2}dx
\le
\frac{\pi}{(KA)^{1/2}}
                             \,  \exp\(
                                       -\dfrac15\l(1+\l)^{-1} K^{2\alpha}
                                    \),
$$
($K$ being one of $mn$, $m$, or $n$; and $\alpha$
being, respectively, $2\eps$, $\eps$, or $\eps$),
one of the $m+n-1$ factors is smaller than $(\pi/Amn)^{1/2}$,
$(\pi/Am)^{1/2}$, or $(\pi/An)^{1/2}$ (respectively, depending
on~$K$) times $\exp\(-\tfrac15\l(1+\l)^{-1}\min(m^{2\eps},n^{2\eps})\)$.
(One may
assure that neither of the indices $j_*$, $k_*$ chosen
for omission in the change of variables coincides with the
index connected to the $h$th inequality of the negation
of definition (\ref{RDef}).)   The factor $\exp(-\tfrac15\cdots)$
is, of course, independent of $h$.  Allowing $N$ values of $\ell$
and $m+n+1$ values of $h$, we have
\begin{align*}
\int_{\U_0\cap\R^c} |F| &\le (m+n+1) N O(\delta) 
n^{O(1)} \,
        \left(\frac{\pi}{Amn}\right)^{1/2}
   \,   \left(\frac{\pi}{Am}\right)^{(n-1)/2}
   \,   \left(\frac{\pi}{An}\right)^{(m-1)/2}  \\
 & ~~~~ \times \exp\(
                        -\dfrac15\l(1+\l)^{-1}\min(m^{2\eps},n^{2\eps})
                   \) \\
            &\le e^{-n^\eps} I_0.
\end{align*}
This inequality
and (\ref{R7}) together imply (\ref{R3}), thus completing
the proof of Theorem~\ref{Jintegral}.
\end{proof}

\nicebreak
\section{Proof of Theorem \ref{thm:main}}\label{section:ProofThm1}

Once again we remind the reader that $F(\bit,\bip)$ denotes the integrand
in equation (\ref{Idef2})
defining $I(m,n)$.  We continue to use
$I_0$ for the quantity defined in equation (\ref{ThmAssumption}).

\medskip

\begin{proof}[Proof of Theorem \ref{thm:main}]

Let $a,b$ be given positive numbers satisfying $a+b<\tfrac12$, and let
$(m,n,s,t)$ be a sequence of 4-tuples such that $m,n\rightarrow\infty$
in such a way that hypothesis~\eqref{Hyp} is satisfied.

We first note that both~\eqref{mainformula2} and
Corollary~\ref{biglambda} are equivalent to~\eqref{mainformula}
for the respective ranges of $\l$ indicated.
Consequently, to prove Theorem~\ref{thm:main} and Corollary~\ref{biglambda}
it suffices to prove~\eqref{mainformula2} for $\lambda=O(n^5)$ and
Corollary~\ref{biglambda} for $\lambda \ge n^5$.

We will divide the proof into two ranges of the parameter~$\l$.
First we assume that $\l=O(n^5)$, where Theorem~\ref{boxing}
can be applied, then we use a different method for
$\l\ge n^5$.

First suppose that $\l=O(n^5)$.
As explained in Section \ref{section:evalint},
\begin{equation}
\label{FiveOne}
 M(m,s;n,t) = \frac{1}{(2\pi)^{m+n}} \(\l^{-\l}(1+\l)^{1+\l}\)^{mn} \,
             \biggl(\,
                      \int_{\R'}F + O(1)\int_{\R^c}|F|
             \biggr),
\end{equation}
where $\R$ is
defined by~\eqref{RDef} and $\R' \supseteq \R$ is defined
by~\eqref{RprimeDef}.
By~\eqref{IvsJ} and Theorem~\ref{Jintegral}, the
first integral on the right side of~\eqref{FiveOne} equals
$2\pi \exp\(\frac12+O(n^{-b})\) I_0$, and by
Theorem~\ref{boxing}, the second equals $O(n^{-1}) I_0$.
These values
yield Theorem~\ref{thm:main} for this case in the
form given by~\eqref{mainformula2}.

\smallskip

For the remainder of the proof, we assume that $\lambda\ge n^5$.
By the first part of this proof, Corollary~\ref{biglambda}
holds for $\l=O(n^5)$; hence,
\begin{equation}\label{ratio1}
  \frac{M(m,2n^6; n,2mn^5)}{M(m,n^6; n,mn^5)}
  = 2^{(m-1)(n-1)}
      \(1+O(n^{-b})\).
\end{equation}
Let $s_0=\lcm(m,n)/m$, $t_0=\lcm(m,n)/n$, and let $\P$ be
the convex polytope of $m\times n$ real, nonnegative matrices
whose rows sum to $s_0$ and whose columns sum to~$t_0$.
For every pair $(s,t)$ such that $ms=nt$ there is an integer
$q$ such that $(s,t)=(qs_0,qt_0)$, and
$M(m,s;n,t)$ equals the number
of integer lattice points in the dilated polytope $q\P$.
The latter count, as a function of the positive integer $q$,
is called the
{\it polytope enumerator} for $\P$, and is denoted $L_{\P}(q)$,
\cite{BeckRobins}.  Thus, we have 
$$
M(m,s;n,t)=L_{\P}(q).
$$
The polytope $\P$, an example of a {\it transportation polytope},
has integral vertices \cite{Brualdi}, and
so $L_{\P}(q)$ is a polynomial in $q$, the {\it Ehrhart 
polynomial}.  The degree of $L_{\P}$ is $d$, the dimension
of the polytope $\P$; in our case, $d=(m{-}1)(n{-}1)$.  By a theorem of
Stanley (an algebraic proof first appeared as Proposition~4.5
in~\cite{stanleyC}; a more geometric proof is given as Theorem~2.1
in~\cite{stanleyB}) the representation of $L_{\P}(q)$ in a
particular basis,
$$
 L_{\P}(q) = \sum_{i=0}^d h_{d-i} \binom{q+i}{d},
$$
has all its coefficients $h_0,\ldots,h_d$ nonnegative.
For $q\ge d$, we have the expansion
\[
 \binom{q+i}{d}
  = \binom{q}{d} \prod_{j=0}^{d-1}\biggl(1+\frac{i}{q-j}\biggr)
  = \binom{q}{d} \negthinspace
     \sum_{X\subseteq[0,d-1]} \frac{i^{|X|}}{\prod_{j\in X} (q-j)},
\]
and so
\[
   L_{\P}(q) = \binom{q}{d} \negthinspace
     \sum_{X\subseteq[0,d-1]} \frac{g(|X|)}{\prod_{j\in X} (q-j)}
\]
where $g(k)=\sum_{i=0}^d h_{d-i}i^k\ge 0$.
Note that $\P$, and therefore $g(k)$ for each $k$, depend only
on $m$ and~$n$.
Also note that $q/\prod_{j\in X} (q-j)$ is decreasing as a
function of~$q$ for $q\ge d$ and any non-empty $X\subseteq[0,d{-}1]$.

Since $q=\l\gcd(m,n)$, we conclude that 
there is a function $\alpha(m,n,\lambda)$ such that
\begin{align}
   M(m,\l n; n,\l m)& =
   \binom{\l \gcd(m,n)}{(m{-}1)(n{-}1)} g(0) 
      \(1 + \alpha(m,n,\l)/\l\) \label{bit1}\\
      \alpha(m,n,\l)&\ge 0
        \text{~~for $q\ge d$}\label{bit2}\\[0.3ex]
   \alpha(m,n,\l)& 
        \text{~~is decreasing in $\l$ for fixed $m,n$ and $q\ge d$}. \label{bit3}
\end{align}
For $\l\ge n^5$, since $d=\tO(n^2)$,
$$
\binom{q}{d} = \frac{q^d}{d!}\(1+\tO(n^{-1})\).
$$
Hence, by~\eqref{bit1}
\begin{equation}\label{ratio2}
  \frac{M(m,2n^6; n,2mn^5)}{M(m,n^6; n,mn^5)}
  = 2^{(m-1)(n-1)} \, \(1+\tO(n^{-1})\) \,
      \frac{1+\tfrac12\alpha(m,n,2n^5)/n^5}{1+\alpha(m,n,n^5)/n^5}.
\end{equation}
Comparing this to~\eqref{ratio1} and noting
from~(\ref{bit2},\,\ref{bit3})
that $0\le \alpha(m,n,2n^5)\le\alpha(m,n,n^5)$,
we conclude that $\alpha(m,n,n^5)=O(n^{5-b})$.
This implies by~\eqref{bit3} that
$\alpha(m,n,\lambda)=O(n^{5-b})$ for $\lambda\ge n^5$.
Using this information about $\alpha(m,n,\l)$ with~\eqref{bit1} gives
\[
   \frac{M(m,\lambda n; n,\lambda m)}{M(m,n^6; n,mn^5)}
= \biggl(\frac{\l}{n^5}\biggr)^{\!(m-1)(n-1)} \(1+O(n^{-b})\).
\]
Combining this with Corollary~\ref{biglambda}
for $\l=n^5$ shows that Corollary~\ref{biglambda} holds for
all $\lambda\ge n^5$.
\end{proof}

\end{document}